\documentstyle[amscd,11pt,leqno]{amsart}

\newtheorem{theorem}{Theorem}[section]
\newtheorem{lemma}[theorem]{Lemma}
\newtheorem{prop}[theorem]{Proposition}
\newtheorem{cor}[theorem]{Corollary}

\newtheorem{question}[theorem]{Question}

\theoremstyle{definition}

\newtheorem{example}[theorem]{Example}

\theoremstyle{remark}
\newtheorem{remark}[theorem]{Remark}

\numberwithin{equation}{section}

\newenvironment{prf}
{\begin{trivlist}\item[]{\bf Proof:\ }} {\qed \end{trivlist}}

\newcommand{\tr}{{\rm tr}}
\renewcommand{\Bbb}{\bf}

\begin{document}

\title[Rigidity of periodic diffeomorphisms]
  {Rigidity of periodic diffeomorphisms on homotopy $K3$ surfaces}


\author{Jin Hong Kim}

\date{
\today; Preliminary version}



\address{Department of Mathematics\\
         KAIST, Kusong-dong, Yusong-gu\\
         Daejon 305--701, South Korea}

\email{jinkim11@@kaist.ac.kr}


\begin{abstract}
In this paper we show that homotopy $K3$ surfaces do not admit a
periodic diffeomorphism of odd prime order $3$ acting trivially on
cohomology. This gives a negative answer for period $3$ to Problem
4.124 in the Kirby's problem list. In addition, we give an
obstruction in terms of the rationality and sign of the spin numbers
to the existence of a periodic diffeomorphism of odd prime order
acting trivially on cohomology of homotopy $K3$ surfaces. The main
strategy is to calculate the Seiberg-Witten invariant for the
trivial spin$^c$ structure in the presence of such a ${\Bbb
Z}_p$-symmetry in two ways: (1) the new interpretation of the
Seiberg-Witten invariants of M. Furuta and F. Fang, and (2) the
theorem of J. Morgan and Z. Szab\' o on the Seiberg-Witten invariant
of homotopy $K3$ surfaces for the trivial Spin$^c$ structure. As a
consequence, we derive a contradiction for any periodic
diffeomorphism of prime order $3$ acting trivially on cohomology of
homotopy $K3$ surfaces.
\end{abstract}

\maketitle

\setlength{\baselineskip}{15pt}

\section{Introduction} 

For compact complex surfaces $X$, it is a well-known question to ask
whether the group $\text{Aut}(X)$ of holomorphic automorphisms of
$X$ acts faithfully on the cohomology ring $H^\ast(X, A)$ with
values in some ring $A$. In general, the answer is negative if the
Lie algebra of $\text{Aut}(X)$ does not reduce to $\{ 0 \}$. For a
concrete example, there exists an Enriques surface $X$ for which
$\text{Aut}(X)$ does not act faithfully on $H^2(X, {\bf Q})$ with
values in the rational numbers. However, for a $K3$ surface $X$ the
automorphism group $\text{Aut}(X)$ is known to act faithfully on the
second cohomology group $H^2(X, {\bf Z})$. In particular, a
holomorphic cyclic ${\Bbb Z}_p$-action of prime order $p$ on a $K3$
surface cannot be (co)homologically trivial. See \cite{Pe} and
references therein for more details.

The aim of this paper is to give a negative answer to the
following Problem 4.124 in the Kirby's problem list \cite{Kir}
which was proposed by Allan Edmonds:

\begin{question} \label{ques1.1}
Do $K3$ surfaces admit a periodic diffeomorphism of prime order
acting trivially on cohomology (or homology)?
\end{question}

The answer to the Question \ref{ques1.1} has already been given
negatively for period $2$ by T. Matumoto \cite{Ma} and independently
by D. Ruberman \cite{Ru}. This shows the rigidity of homologically
trivial actions of prime order on $K3$ surfaces.

In the present paper we attempt to answer the Question \ref{ques1.1}
even on \emph{homotopy} $K3$ surfaces, not just $K3$ surfaces, by a
contradiction. In a little more detail, we first assume that the
spin number of a periodic diffeomorphism $\tau$ of odd prime order
satisfying $b_+(X/\tau)=3$ is both rational and negative. Then we
show in Theorem \ref{thm4.1} that the Seiberg-Witten invariant for
the trivial spin$^c$ structure vanishes identically. To show such a
vanishing result, we use the new $K$-theoretic interpretation of the
Seiberg-Witten invariants introduced by Furuta and nicely developed
by Fang. On the other hand, Morgan and Szab\' o have already proved
in \cite{M-S} that the Seiberg-Witten invariant on homotopy $K3$
surfaces for the trivial spin$^c$ structure should be $\pm 1$ mod 2.
This implies an obvious contradiction for odd prime period under the
assumption that the spin number is both negative and rational.

In case of prime period $3$, we can show further that the spin
number is indeed rational (see Corollary \ref{cor3.2}). Moreover, it
is also true that a homotopy $K3$ surface admits a periodic
diffeomorphism of order $3$ which acts trivially on cohomology only
if its spin number is negative. Thus we can give a negative answer
to Question \ref{ques1.1} as follows.

\begin{theorem} \label{thm1.1}
Let $X$ be a homotopy $K3$ surface, and let $\tau: X\to X$ be a
periodic diffeomorphism of order $3$. Then $\tau$ cannot act
trivially on cohomology.
\end{theorem}

Here we remark that the condition that the periodic diffeomorphism
acts trivially on cohomology cannot be replaced by the condition
$b_+(X/\tau)=3$. Indeed, there exists an action of ${\bf Z}_3$ on a
certain quartic $K3$ surface $X$ in ${\bf CP}^3$ which acts
trivially on $H^{2}_+(X; {\bf R})$, as S. Mukai shows on page 187 of
\cite{Mu}. To be precise, there exists an action of $PSL(2, {\bf
F}_7)$ of a certain $K3$ surface in ${\bf CP}^3$. Since the group
$PSL(2, {\bf F}_7)$ is simple and the subgroup of $PSL(2, {\bf
F}_7)$ which acts trivially on $H^0(X, \Omega_X^2)$ is normal, the
action on $H^0(X, \Omega_X^2)$ is actually trivial. (See also
\cite{Ni}.) Moreover, since the the action of $PSL(2, {\bf F}_7)$
preserves the K\" ahler form of ${\bf CP}^3$, it preserves that of
the $K3$ surface. Thus the action on $H^2_+(X)\otimes {\bf C}$ is
trivial. On the other hand, since $PSL(2, {\bf F}_7)$ contains a
copy of ${\bf Z}_3$, there exists such a ${\bf Z}_3$ action on a
$K3$ surface.

It will be clear from Section 4 that we can give an obstruction to
the existence of a periodic diffeomorphism of odd prime order
acting trivially on cohomology of homotopy $K3$ surfaces as
follows.

\begin{theorem} \label{thm1.2}
Let $X$ be a homotopy $K3$ surface, and let $\tau: X\to X$ be a
periodic diffeomorphism of odd prime order $p$. Assume that the spin
number $\text{\rm Spin}(\hat\tau, X)$ is both rational and negative.
Then $\tau$ cannot act trivially on the self-dual part $H^{2}_+(X;
{\bf R})$ of the second cohomology group.
\end{theorem}

As an immediate corollary, we can give the following

\begin{cor}
Let $X$ be a homotopy $K3$ surface, and let $\tau: X\to X$ be a
periodic diffeomorphism of odd prime order $p$. Assume that the
spin number $\text{\rm Spin}(\hat\tau, X)$ is rational and
negative. Then $\tau$ cannot act trivially on cohomology.
\end{cor}

For a detailed definition of the spin number $\text{\rm
Spin}(\hat\tau, X)$, see Section 3. In contrast to the above
negative result, it is interesting to note that Edmonds \cite{Ed}
has constructed locally linear homologically trivial ${\Bbb
Z}_p$-actions of odd prime order $p$ on any closed simply
connected \emph{topological} 4-manifold.

Recently Chen and Kwasik announced, among many other things, a
negative result in \cite{C-K} to Question \ref{ques1.1} under the
stronger assumption that a \emph{symplectic} action of odd prime
order on $K3$ surfaces acts trivially on cohomology. Their method is
quite different from ours. In fact, they use a pseudo-holomorphic
curve theory for symplectic 4-orbifolds in order to obtain an
equivariant version of the Taubes' well-known theorem on symplectic
4-manifolds of \cite{Ta}. Moreover, the claim that the fixed points
of such symplectic actions on $K3$ surfaces like holomorphic actions
are always isolated is a crucial ingredient in their proof. In
general, the fixed-point set of smooth actions of odd prime order on
homotopy $K3$ surfaces may or may not be isolated, 2-dimensional or
a union of isolated points and 2-dimensional submanifolds. This
greatly complicates the proof of Theorem \ref{thm1.2} as the present
paper shows. (See \cite{C-K}, \cite{Ch} and \cite{Ta} for more
details.)

We organize this paper as follows. In Section 2, we set up basic
notations and collect some important theorems necessary for the
later sections. In Section 3, we recall the general formula for the
spin numbers, and then we see that the spin number is always real.
In the same section we show a very important relationship between
the differences of the dimensions of the eigenspaces under the
action of ${\Bbb Z}_p$ on the representation spaces induced from the
spinor vector bundles. Section 4 is devoted to the proof of Theorem
\ref{thm4.1} which yields immediately the proof of Theorem
\ref{thm1.2}. In the same Section 4 we give a complete proof of
Theorem \ref{thm1.2}. Finally in Section 5, we give a proof of
Theorem \ref{thm1.1}.

\smallskip\smallskip

\noindent{\bf Acknowledgements:} The author is grateful to Professor
Ian Hambleton for telling the work of W. Chen and S. Kwasik to him
by email and Professor Weimin Chen for sending his manuscript
\cite{C-K}.

\smallskip\smallskip

\section{Preliminaries}

In this section, we briefly review the Seiberg-Witten equations
and the new interpretation of Seiberg-Witten invariants by Furuta
and Fang. (See \cite{Fu1, Fu, Fa, Br, Kim1} for more details.)

\subsection{Seiberg-Witten invariants:}
Let $c$ be a spin$^c$ structure on $X$ whose associated line
bundle is $L$, and $S_{\pm}$ denote the positive or negative
spinor vector bundles. Let ${\cal A}$ be the space of connections
on $L$, and let $\Omega^+(X)$ be the space of self-dual 2-forms.
We denote by ${\cal G}$ the group of gauge transformations on $L$.
Then using the Seiberg-Witten equations we have a ${\cal
G}$-equivariant map
\[
f: {\cal A}\times \Gamma(S_+)\to i \Omega^+(X)\times \Gamma(S_-),
\quad (A, \phi)\mapsto (F_A^+ -q(\phi), D_A\phi),
\]
where $q(\phi)=\phi\otimes \phi^\ast-\frac{|\phi|^2}{2}\text{Id}$.
For the sake of convenience, we assume that $\Gamma(S_+)$ and
${\cal A}$ are completed with a Sobolev norm $L_2^4$ and that
$\Gamma(S_-)$ and $\Omega^+(X)$ are completed with a Sobolev norm
$L_2^3$.

We fix a connection $A_0$. Since $A=A_0+a$ and the stabilizer of
$A_0$ is $S^1$, we have an $S^1$-equivariant map
\[
f: i\Omega_c^1(X)\times \Gamma(S_+)\to i\Omega^+(X)\times
\Gamma(S_-), \quad (a, \phi)\mapsto (d^+a-q(\phi),
D_{A_0}\phi+a\cdot \phi),
\]
where $\Omega_c^1(X)$ is the  subspace of 1-forms in the
$\text{\rm Ker}\, d^\ast$.

To construct the finite dimensional approximation of Furuta, we let
$U_0=\Gamma(S_+)$ and $U'_0=\Gamma(S_-)$, and let
$V_0=i\Omega_c^1(X)$ and $V'_0=i\Omega^+(X)$. For each positive real
number $\lambda$, we define $U_\lambda$ (resp. $U'_\lambda$) the
vector space spanned by eigenvectors of the  operator $D_{A_0}^\ast
D_{A_0}$ (resp. $D_{A_0} D_{A_0}^\ast$) with eiegenvalues less than
or equal to $\lambda$.  Similarly we define $V_\lambda$ and
$V_\lambda'$ using the elliptic operator $d^+$.

Let $p_\lambda$ denote the $L_2$-orthogonal projection of
$U_0\oplus V_0$ onto $U'_\lambda\oplus V'_\lambda$.  Using the
restriction of $f$ and the projection $p_\lambda$, we have an
$S^1$-equivariant map
\[
f_\lambda: U_\lambda\oplus V_\lambda\to U'_\lambda\oplus
V'_\lambda,
\]
where $S^1$ acts on $\Omega^i(X)$ ($i=1,2$) trivially and on
$\Gamma(S_\pm)$ by the complex multiplication.

Let $W_\lambda'=U_\lambda'\oplus V'_\lambda$ and $\nu_0=-F_{A_0}^+
+\mu$ for a generic 2-form $\mu$. Using the compactness of the
Seiberg-Witten moduli space, Furuta showed the following lemma (see
also Lemma 2.1-2 in \cite{Fa}):

\begin{lemma}[\cite{Fu}] \label{lem2.1}
For a generic parameter $\mu$ as above and a sufficiently large
positive $R$, there exists a positive number $\Lambda$ such that
$f_\lambda^{-1}(\nu_0)$ and $(0\times V_\lambda)\cap B_\lambda(R)$
do not intersect for $\lambda\ge \Lambda$, where $B_\lambda(R)$ is
the ball with radius $R$ in $W_\lambda$.
\end{lemma}

From now on, we assume that $\lambda$ is greater than or equal to
$\Lambda$. Then Lemma \ref{lem2.1} gives rise to an
$S^1$-equivariant map
\[
f_\lambda: (B_\lambda, \partial B_\lambda\cup ((0\times
V_\lambda)\cap B_\lambda))\to (W'_\lambda, W'_\lambda -
B_\epsilon(\nu_0)),
\]
where $B$'s denote balls. Now, passing to the quotients of the
previous map $f_\lambda$, we get an $S^1$-equivariant map
\begin{equation*}
f_\lambda: S^{V_\lambda\oplus {\Bbb R}}\wedge S(U_\lambda)\to
S^{W'_\lambda},
\end{equation*}
where $S^{V_\lambda\oplus {\Bbb R}}$ and $S^{W'_\lambda}$ are
called the Thom spaces.

Let $m=\dim_{\Bbb C} U_\lambda, m'=\dim_{\Bbb R} V_\lambda,
n=\dim_{\Bbb C} U'_\lambda$, and $n'=\dim_{\Bbb R} V'_\lambda$.
Then we have
\[
\frac{1}{4}(c_1(L)^2-\sigma(X))=2m-2n \quad\text{and}\quad
-b_2^+(X)=m'-n'.
\]
Thus we have
$2d=\frac{1}{4}(c_1(L)^2-(2\chi(X)+3\sigma(X)))=2m+m'-(2n+n'+1)$.
By the Thom isomorphism theorem we have the following commutative
diagram
\begin{equation*}
\begin{CD}
H^{2n+n'}_{S^1}(S^{W'_\lambda}, {\Bbb Z}) @> f^\ast_\lambda>>
H^{2n+n'}_{S^1}(S^{V_\lambda\oplus {\Bbb R}}\wedge S(U_\lambda),
{\Bbb Z})\\
@.      @VV \tau V\\
@.  H^{2(m-1-d)}_{S^1}(S(U_\lambda), {\Bbb Z}),\\
\end{CD}
\end{equation*}
where $\tau$ is the suspension isomorphism and all the cohomology
appearing in the diagram means $S^1$-equivariant cohomology.

Let $\Phi\in H^{2n+n'}_{S^1}(S^{W'_\lambda}, {\Bbb Z})$ be the
equivariant Thom class. Then from the above diagram we have an
element $\theta_{f_\lambda}\in H^{2(m-1-d)}_{S^1}(S(U_\lambda),
{\Bbb Z})$ such that
\[
f_\lambda^\ast(\Phi)=\tau^{-1}(\theta_{f_\lambda}).
\]
Since $S^1$ acts freely on $S(U_\lambda)$, we see that
\[
H^{2(m-1-d)}_{S^1}(S(U_\lambda), {\Bbb Z})=H^{2(m-1-d)} ({\Bbb
CP}^{m-1}, {\Bbb Z}).
\]
Note also that since $H^{\ast} ({\Bbb CP}^{m-1}, {\Bbb Z})={\Bbb
Z}[x]/x^m$, we can write $\theta_{f_\lambda}=a_{X, c} x^{m-1-d}$
for some $a_{X, c}\in {\Bbb Z}$.

A cyclic ${\Bbb Z}_q$-action on $X$ ($q$ is not necessarily prime)
is called a \emph{spin$^c$ action} (or \emph{preserves a spin$^c$
structure}) if the generator of the action $g: X\to X$ lifts to a
spin$^c$ bundle $\hat g: P_{Spin^c}\to P_{Spin^c}$.  Such an
action is of even type if $\hat g$ has order $q$ and is of odd
type if $\hat g$ has order $2q$. In particular, if the spin$^c$
bundle is a spin bundle, the action is called a \emph{spin
action}.  According to \cite{Fa}, a cyclic ${\Bbb Z}_p$-action of
odd prime $p$ is of even type if and only if the associated line
bundle to the spin$^c$ structure is a ${\Bbb Z}_p$-$U(1)$-bundle.
Since the line bundle associated to the trivial spin$^c$ structure
is always a ${\Bbb Z}_p$-$U(1)$-bundle, any cyclic spin action of
odd prime order  is of even type.

From now on, we denote by $\hat{\Bbb Z}_q$ the group generated by
$\hat g$. Thus if $X$ has an action ${\Bbb Z}_q$ preserving the
spin$^c$ structure $c$, then we have an $S^1\times \hat {\Bbb
Z}_q$-equivariant map $f=f_\lambda: S^{V_\lambda\oplus {\Bbb
R}}\wedge S(U_\lambda)\to S^{W'_\lambda}$. Moreover, by suspending
$f$ by $S^{V_\lambda\oplus {\Bbb R}}$ we obtain a map, still
denoted $f$,
\begin{equation} \label{eq2.1}
f: S^{(V_\lambda\oplus {\Bbb R})\otimes {\Bbb C}}\wedge
S(U_\lambda)\to S^{W'_\lambda\oplus V_\lambda\oplus {\Bbb R}}.
\end{equation}
Applying the $K_{S^1\times \hat{\Bbb Z}_q}$-theory we get a
$\beta_f\in K_{S^\times \hat{\Bbb Z}_q}(S(U_\lambda))$ such that
\[
f^\ast(\tau_{W'_\lambda\oplus V_\lambda\oplus {\Bbb R}})=\beta_f
\tau_{(V_\lambda\oplus {\Bbb R})\otimes {\Bbb C}},
\]
where $\tau_{W'_\lambda\oplus V_\lambda\oplus {\Bbb R}}$ and
$\tau_{(V_\lambda\oplus {\Bbb R})\otimes {\Bbb C}}$ are the
$K$-theory Thom classes. Then we can summarize the results of
Furuta and Fang in the following theorem (see Theorems 2.3 and 2.4
and Proposition 4.3 in \cite{Fa}).

\begin{theorem}[Furuta, Fang] \label{thm2.1}
Let $X$ be a smooth 4-manifold with $b_1(X)=0$ and that
$b_2^+(X)\ge 2$. Let $c$ denote spin$^c$ structure on $X$.
\begin{itemize}
\item[(1)]
For a sufficiently large $\lambda\ge \Lambda$, the Seiberg-Witten
invariants satisfy $SW(X, c)=a_{X, c}$. Furthermore, if $X$ has an
action ${\Bbb Z}_q$ preserving the spin$^c$ structure $c$ and
$H^2_+(X/{\Bbb Z}_q, {\Bbb R})\ne 0$, then there exists an
$S^1\times \hat{\Bbb Z}_q$-equivariant map
\[
f_\lambda: S^{V_{\lambda}\oplus {\Bbb R}}\wedge S(U_\lambda)\to
S^{W_{\lambda}'}
\]
and $SW(X, c)=a_{X, c}$.

\item[(2)] Let $t$ be the standard 1-dimensional complex representation
of $S^1$, and let $T=1-t$. Then $\beta=\beta_f$ satisfies the
following identity:
\begin{equation*}
\beta_f(t)=(-1)^n a_{X, c} \left(
\frac{\log(1+T)}{T}\right)^{\frac{1}{2}(b_+(M)-1)} T^{m-d-1},\
T^{m-d}=0.
\end{equation*}
In particular, we have $\beta=\pm SW(X, c)T^{m-d-1}$.
\end{itemize}
\end{theorem}

\subsection{The tom Dieck's character formula:}
We also need to use the tom Dieck's formula in Section 3. we
explain it briefly. Since any spin cyclic actions of odd prime
order are of even type, we consider spin actions only of even
type.

Recall first that the group $Pin(2)$ has one non-trivial one
dimensional representation $\tilde 1$ and  has a countable series
of 2-dimensional irreducible representations $h_1, h_2,\ldots$. In
particular, the representation $h_1=h$ is the restriction of the
standard representation of $SU(2)$ to $Pin(2)\subset SU(2)$. The
representation ring $R({\Bbb Z}_p)$ is isomorphic to the group
ring ${\Bbb Z}({\Bbb Z}_p)$ which is generated by the standard
one-dimensional representation $\xi$. Thus as a ${\Bbb Z}$-module,
$R({\Bbb Z}_p)$ is generated by $1,\xi,\ldots, \xi^{p-1}$.

For the sake of simplicity, let $V=(U_\Lambda\oplus
V_\Lambda)\otimes_{\bf R} {\Bbb C}$, $W=(U_\Lambda'\oplus
V_\Lambda')\otimes_{\bf R} {\Bbb C}$, and $G=Pin(2)\times {\Bbb
Z}_p$. In the presence of a spin action ${\Bbb Z}_p$ of odd prime
order $p$, the complex index of $D$ is given by
\begin{equation*}
[V]-[W]=k(\xi)h-t(\xi)\tilde 1\in R(G),
\end{equation*}
where $k(\xi)=k_0+k_1\xi+\ldots+k_{p-1}\xi^{p-1}$,
$t(\xi)=t_0+t_1\xi+\ldots+t_{p-1}\xi^{p-1}$ with the properties
$t_0+t_1+\ldots+t_{p-1}=b_+(X)$ and
$k_0+k_1+\ldots+k_{p-1}=-\frac{\sigma(X)}{8}$.

Let $BV$ and $BW$ denote balls in $V$ and $W$. Then it follows
from Lemma \ref{lem2.1} that there exists a $G$-equivariant map
$f$ preserving the boundaries $SV$ and $SW$ of $BV$ and $BW$,
respectively:
\begin{equation*}
f=f_\Lambda: (BV, SV)\to (BW, SW).
\end{equation*}
Let $K_G(V)$ denote $K_G(BV, SV)$. Similarly define $K_G(W)$. Then
$K_G(V)$ (resp. $K_G(W)$) is a free $R(G)$-modules with one
generator $\lambda(V)$ (resp. $\lambda(W)$), called the \emph{Bott
class}. Now, applying $K$-theory functor to $f$ we get a map
\[
f^\ast: K_G(W)\to K_G(V)
\]
with a unique element $\alpha_f$, called the \emph{$K$-theoretic
degree of $f$}, satisfying the equation
\begin{equation} \label{eq2.2}
f^\ast(\lambda(W))=\alpha_f\cdot \lambda(V).
\end{equation}
Let $V_g$ and $W_g$ denote the subspaces of $V$ and $W$ fixed by
an element $g\in G$, and $V_g^\perp$ and $W_g^\perp$ denote the
their orthogonal complements. Let $f^g: V_g\to W_g$ be the
restriction of $f$ to $V_g$, and let $d(f^g)$ denote the
topological degree of $f^g$. Then the tom Dieck's character
formula says that we have
\begin{equation*}
\tr_g(\alpha_f)=d(f^g)\tr_g\left(\sum_{i=0}^\infty (-1)^i
\Lambda^i(W_g^\perp-V_g^\perp)\right),
\end{equation*}
where $\tr_g$ is the trace of the action of an element $g\in G$.
Note also that the topological degree $d(f^g)$ is by definition
zero, if $\dim(V_g)\ne \dim(W_g)$.

\section{Cyclic Group Actions and Spin Numbers} 

In this section, we prove the important lemmas about the spin
numbers which are essential to prove our main result.

As before, let $X$ denote a homotopy $K3$ surface and let $\tau$ be
a periodic diffeomorphism of odd prime order $p$, unless stated
otherwise. Let $\sigma(X)$ denote the signature of $X$. As remarked
in Subsection 2.1, it suffices to consider spin actions $\tau$ of
even type.

First recall that the spin number for the lifting $\hat \tau$  is
defined to be
\[
\text{\rm Spin}(\hat\tau, X)=\text{ind}_{\hat\tau}
D=\tr(\hat\tau|_{\text{ker}D})-\tr(\hat\tau|_{\text{coker}D}),
\]
where $D$ denotes the Dirac operator as before. These spin numbers
can be calculated in terms of the fixed point set $X^\tau$ by the
general Lefschetz formula (Theorem 3.9 in \cite{A-S}).

In more detail, for each $x\in X^\tau$, the fiber $N_x^\tau$ of the
normal bundle $N^\tau$ is a real ${\bf Z}_p$-module. Since ${\bf
Z}_p$ is cyclic of odd prime order, its real irreducible
representation is of the form
\begin{equation*}
\tau\mapsto
\begin{pmatrix}
\cos\theta & -\sin\theta\\
\sin\theta& \cos\theta
\end{pmatrix}.
\end{equation*}
It is important to note that the above representation given by
$\theta$ and $-\theta$ are equivalent. Thus the fiber $N_x^\tau$ of
the normal bundle has a canonical decomposition
\[
N_x^\tau= \sum_\mu N_x^\tau(\mu),
\]
where $\mu$'s are the complex numbers of absolute value $1$ with
positive imaginary part and where $N_x^\tau(\mu)$ has a complex
structure in which $\tau$ acts by multiplication with $\mu$. Then we
have the following Lefschetz theorem about the spin numbers.

\begin{theorem} \label{thm3.1}
Let $\tau: X\to X$ be an isometry of the homotopy $K3$ surface $X$
of odd prime order $p$. Let the fixed-point set $X^\tau$ consist of
isolated points $\{ P_j \}$ and connected 2-manifolds $\{ F_k\}$.
For each $j$, let the action of $\tau$ on the tangent space at $P_j$
be given by the matrix
\begin{equation*}
\begin{pmatrix}
\cos\alpha_j & -\sin\alpha_j\\
\sin\alpha_j & \cos\alpha_j
\end{pmatrix}
\oplus
\begin{pmatrix}
\cos\beta_j & -\sin\beta_j\\
\sin\beta_j & \cos\beta_j
\end{pmatrix}
\end{equation*}
relative to an oriented basis, where $\alpha_j$ and $\beta_j$ denote
${2\pi l_{\alpha_j}}/{p}$ and ${2\pi l_{\beta_j}}/p$ $(0< \alpha_j,
\beta_j<\pi)$, respectively. For each $k$, let $\tau$ act on the
normal bundle $N^\tau_k$ of $Y_k$ by multiplication with
$e^{i\theta_k}$, where $\theta_k={2 \pi l_{\theta_k}}/p$ $(0<
\theta_k <\pi)$. Let $\hat\tau$ be a lifting of $\tau$ preserving
the trivial spin$^c$ structure. Then we have the following formula
for spin number:
\begin{equation} \label{eq3.1}
\begin{split}
\text{\rm Spin}(\hat\tau, X)&=-\frac{1}{4}\sum_{P_j} \epsilon(P_j,
\hat\tau) \csc(\alpha_j/2)\csc(\beta_j/2)\\
 &+\frac{1}{4}\sum_{F_k} \epsilon(F_k, \hat\tau)\cos(\theta_k/2)
 \csc^2(\theta_k/2)
\langle [F_k], [F_k]\rangle,
\end{split}
\end{equation}
where $\epsilon(P_j, \hat\tau)$ and $\epsilon(F_k, \hat\tau)$ are
$\pm 1$, depending on the action of $\tau$ on the spin bundle.
\end{theorem}

The proof of the above theorem follows from the Atiyah-Singer
$G$-spin theorem in the literature (e.g., see \cite{Hi},
\cite{AH}, Theorem 8.35 in \cite{A-B} and Theorem 14.11 in
\cite{L-M}).

Then as a corollary we can show the following result.

\begin{cor} \label{cor3.1}
Let $X$ be a homotopy $K3$ surface, and let $\tau$ be a periodic
diffeomorphism of odd prime order $p$ on $X$. Let $\hat\tau$ be a
lifting of $\tau$ preserving the trivial spin$^c$ structure. Then
the following holds:
\begin{itemize}
\item[\rm (1)] The spin number $\text{Spin}(\hat\tau, X)$ is
always real.

\item[\rm (2)] For each $j=1,2,\ldots, p-1$,
two spin numbers $\text{\rm Spin}(\hat\tau^j, X)$ and $\text{\rm
Spin}(\hat\tau^{p-j}, X)$ equal to each other.
\end{itemize}
\end{cor}

\begin{prf}
The proof (a) is immediate from Theorem \ref{thm3.1}.

For the proof of (b), as remarked earlier it suffices to note that
the representation of the normal bundle of a fixed-point set given
by $\theta$ and $-\theta$ are equivalent.
\end{prf}

As another immediate corollary, we have the following result. This
will be useful to prove Theorem \ref{thm1.1} which is the case for
prime order 3.

\begin{cor} \label{cor3.2}
Let $X$ be a homotopy $K3$ surface, and let $\tau$ be a periodic
diffeomorphism of odd prime order $3$ on $X$. Let $\hat\tau$ be a
lifting of $\tau$ preserving the trivial spin$^c$ structure. Then
the spin number $\text{\rm Spin}(\hat\tau, X)$ is always rational.
\end{cor}

\begin{prf}
If the prime order $p$ equals 3 then all of $\alpha_j$, $\beta_j$,
and $\theta_k$ are just $2\pi/3$. But then
$\csc(\alpha_j/2)\csc(\beta_j/2)$ is just $\frac{4}{3}$ and
$\cos(\theta_k/2)\csc^2(\theta_k/2)$ is $\frac23$. Hence it follows
from \eqref{eq3.1} that the spin number $\text{Spin}(\hat\tau, X)$
should be
\begin{equation*}
\text{\rm Spin}(\hat\tau, X)=\frac{1}{6}\sum_{F_k} \pm \langle
[F_k], [F_k]\rangle+\frac13\sum_{P_j}\pm 1 ,
\end{equation*}
which is rational, as required.
\end{prf}

Next, let $m_i$ ($i=0,1,2,\ldots, p-1$) denote the dimensions of the
$\nu^i$-eigenspaces of a generator of the ${\Bbb Z}_p$-action on
$U_\Lambda$, where $\nu$ is a generator of ${\Bbb Z}_p$. We also
define $n_i$ in a similar way for $U_\Lambda'$. Thus $k_i=m_i-n_i$
for $i=0,1,2,\ldots, p-1$. Then we need the following series of
lemmas:

\begin{lemma} \label{lem3.1}
Let $X$ be a homotopy $K3$ surface, and let $X$ admit a periodic
diffeomorphism $\tau$ of odd prime order $p$ satisfying
$b_+(X/\tau)=3$. Then we have $k_0 \le 2$.
\end{lemma}

\begin{prf}
Since the difference of two spin structures is naturally an element
of $H^1(X, {\Bbb Z}_2)=0$, $\tau$ clearly generates a spin action of
odd prime order $p$ that is of even type.

From the previous section, recall
\[
[V]-[W]=k(\xi) h- t(\xi)\tilde 1,
\]
where $k(\xi)=k_0+k_1\xi+\ldots+k_{p-1}\xi^{p-1}$,
$t(\xi)=t_0+t_1\xi+\ldots+t_{p-1}\xi^{p-1}$ with the properties
\begin{equation*}
\begin{split}
t_0+t_1+\ldots+t_{p-1} &=b_+(X)=3,\\
k_0+k_1+\ldots+k_{p-1} &=-\frac{\sigma(X)}{8}=2.
\end{split}
\end{equation*}
Since $b_+(X)=b_+(X/\tau)$ and $\tau$ is of odd prime order, we
should have $t_0=b_+(X/\tau)=3$, and so $t_1+\ldots+t_{p-1}=0$. Let
$\alpha=\alpha_f$ be the $K$-theoretic degree of $f$ in
\eqref{eq2.2} of the form
\[
\alpha_f=\alpha_0(\xi)+\tilde\alpha_0(\xi)\tilde 1
+\sum_{i=1}^\infty \alpha_i(\xi) h_i
\]
Now, we want to compute $\alpha=\alpha_f$. To do so, note first
that $\theta$ and $\theta\nu$ act non-trivially on $h$ and
trivially on $\tilde 1$, where $\theta\in S^1$ is an element
generating a dense subgroup of $S^1$ and $\nu\in {\Bbb Z}_p$ is a
generator as before. Thus we have
\begin{equation*}
\begin{split}
\dim V_\theta - \dim W_\theta&=-(t_0+t_1+\ldots+t_{p-1})=-b_+(X)<0\\
\dim V_{\theta\nu} - \dim
W_{\theta\nu}&=-t_0=-b_+(X/\tau)=-b_+(X)<0.
\end{split}
\end{equation*}
Hence the topological degrees $d(f^\theta)=d(f^{\theta\nu})=0$, so
we obtain $\tr_\theta(\alpha)=\tr_{\theta\nu}(\alpha)=0$. This
implies
\begin{equation*}
\alpha_0(\nu)+\tilde\alpha_0(\nu)=0, \quad \alpha_i(\nu)=0\quad
i\ge 1.
\end{equation*}
Therefore $\tilde\alpha_0=-\alpha_0$ and $\alpha_i=0$ for all
$i\ge 1$, i.e., $\alpha=\alpha_0(\xi)(1-\tilde 1)$.

Set $\alpha_0(\xi)=a_0+a_1\xi+\ldots+a_{p-1}\xi^{p-1}$. Since $J\nu$
acts non-trivially on $\xi \tilde 1, \xi h$, and $h$, we have $\dim
V_{J\nu}-\dim W_{J\nu}=0$. Thus, by definition, $d(f^{J\nu})=1$.
Note also that
\begin{equation} \label{eq3.2}
\tr_{J\nu}(\alpha)=\tr_{J\nu}\left(\alpha_0(\xi)(1-\tilde
1)\right)=2(a_0+a_1\nu+\ldots+a_{p-1}\nu^{p-1}).
\end{equation}
On the other hand, for each $j=1, 2, \ldots, p-1$ we get
\begin{equation} \label{eq3.3}
\begin{split}
&\tr_{J\nu^j}\left(\sum_{i=0}^\infty (-1)^i
\Lambda^i((t_0+t_1\xi+\ldots+t_{p-1}\xi^{p-1})\tilde 1\right.\\
&\left.\qquad -(k_0+k_1\xi+\ldots+k_{p-1}\xi^{p-1})h)\right)\\
&=\tr_{J\nu^j}\left((1-\tilde 1)^{t_0}(1-\xi^j\tilde 1)^{t_1}\ldots
(1-\xi^{j(p-1)}\tilde 1)^{t_{p-1}}\right.\\
&\left.\cdot (2-h)^{-k_0}(1+\xi^{2j}-\xi h)^{-k_1}\ldots
(1+\xi^{j(p-2)}-\xi^{j(p-1)}h)^{-k_{p-1}}\right)\\
&=2^{3-k_0}(1+\nu^j)^{t_1}(1+\nu^{2j})^{t_2}\ldots (1+\nu^{j(p-1)})^{t_{p-1}}\\
&\cdot (1+\nu^{2j})^{-k_1}(1+\nu^{3j})^{-k_2}\ldots (1+\nu^{p-2j})^{-k_{p-1}}.\\
\end{split}
\end{equation}
where we used $\tr_J(\tilde 1)=-1$, and $\tr_J(h)=0$ in the second
equality. Hence it follows from \eqref{eq3.2} and \eqref{eq3.3} that
$k_0\le 2$. Indeed, if we multiply $p-1$ equations from
\eqref{eq3.1}, we obtain
\begin{equation*}
\begin{split}
&2^{(p-1)(3-k_0)} \left(
\prod_{j=1}^{p-1}(1+\nu^j)\right)^{t_1}\ldots \left(
\prod_{j=1}^{p-1}(1+\nu^{p-2j})\right)^{-k_{p-1}}\\
&= 2^{p-1}\prod_{j=1}^{p-1}\left(\sum_{i=1}^{p-1} a_i
\nu^{ji}\right).
\end{split}
\end{equation*}
Since $\prod_{j=1}^{p-1}(1+\nu^j)=1$, the above equation should be
of the form
\begin{equation*}
\begin{split}
2^{(p-1)(2-k_0)} &= \prod_{j=1}^{p-1}\left( \sum_{i=1}^{p-1} a_i \nu^{ji}\right)\\
&= c_0+c_1\nu+\ldots + c_{p-1} \nu^{p-1},
\end{split}
\end{equation*}
where $c_0, c_1, \ldots, c_{p-1}$ are some integers. Thus, if
$k_0\ge 3$ then we have
\[
1= 2^{(p-1)(k_0-2)}(c_0+c_1\nu+\ldots + c_{p-1} \nu^{p-1})\equiv 0,
\quad\text{mod}\ 2,
\]
which is a contradiction. This completes the proof.
\end{prf}

\begin{remark} \label{rmk3.1}
Let $X$ be a homotopy $K3$ surface, and $X$ admit a periodic
diffeomorphism $\tau$ of odd prime order $p$ satisfying
$b_+(X/\tau)=3$ as in Lemma \ref{lem3.1}. Then we can also use the
lifting $e^{\frac{2\pi qi}{p}} \hat \tau$ instead of $\hat\tau$,
where $0\le q < p$. Let $m_i^q$ ($i=0,1,2,\ldots, p-1$) denote the
dimensions of the $\nu^{i+q}$-eigenspaces of a generator of the
${\Bbb Z}_p$-action on $U_\Lambda$, where $\nu$ is a generator of
${\Bbb Z}_p$. We also define $n_i^q$ in a similar way for
$U_\Lambda'$. Thus we have $k_i^q=m_i^q-n_i^q$ for $i=0,1,2,\ldots,
p-1$. Then applying the same arguments as in Lemma \ref{lem3.1}
implies $k_0^q \le 2$ for each $0 \le q < p$. Since $k_i^q=k_{i+q}$
(here the subscripts are labelled modulus $p$), we can conclude that
$k_i \le 2$ for all $0 \le i \le p-1$.
\end{remark}

If the spin number is rational and non-negative, we can describe
$k_i$ more precisely as follows.

\begin{lemma} \label{lem3.2}
Let $X$ be a homotopy $K3$ surface, and let $X$ admit a periodic
diffeomorphism $\tau$ of odd prime order $p$ satisfying
$b_+(X/\tau)=3$. Assume that the spin number $\text{\rm
Spin}(\hat\tau, X)$ is both rational and non-negative. Then we have
\[
k_0=2\ \text{and}\ k_1=k_2=\ldots=k_{p-1}=0.
\]
\end{lemma}

\begin{prf}
We continue to use the same notations as in the proof of Lemma
\ref{lem3.1}.

To show it, note first that we have
\begin{equation} \label{eq3.4}
\text{Spin}(\hat\tau,X) =k_0+k_1\nu+\ldots + k_{p-1}\nu^{p-1}.
\end{equation}
Now, if we use the formula \eqref{eq3.4} for the spin number
$\text{Spin}(\hat\tau, X)$ and the relation
$1+\nu+\nu^2+\ldots+\nu^{p-1}=0$, we obtain
\begin{equation} \label{eq3.5}
\begin{split}
\text{Spin}(\hat\tau,X) &=k_0+k_1\nu+\ldots +
k_{p-1}\nu^{p-1}\\
&=k_0+k_1\nu+\ldots+k_{p-2}\nu^{p-2}-k_{p-1}(1+\nu+\ldots+\nu^{p-2})\\
&=k_0-k_{p-1}+(k_1-k_{p-1})\nu+\ldots+ (k_{p-2}-k_{p-1})\nu^{p-2}.
\end{split}
\end{equation}
Since the spin number $\text{Spin}(\hat\tau, X)$ is rational by
assumption, the equation \eqref{eq3.5} has rational coefficients.
Thus all the coefficients of \eqref{eq3.5} should vanish, since the
polynomial $x^{p-1}+x^{p-2}+\ldots+x+1$ is irreducible over rational
numbers. In particular, we have
\[
k_1=k_2=\ldots=k_{p-1},
\]
as required.

Note also from the equation \eqref{eq3.4} and the relation
$1+\nu+\ldots+\nu^{p-1}=0$ that we have
\begin{equation*}
\begin{split}
(p-1)k_0 &=(p-1)k_1+(p-1)\cdot\text{Spin}(\hat\tau, X)\\
&=2-k_0+ (p-1)\cdot\text{Spin}(\hat\tau, X).
\end{split}
\end{equation*}
Thus the spin number satisfies
\begin{equation} \label{eq3.8}
0\le\text{Spin}(\hat\tau, X)= \frac{p k_0}{(p-1)}-\frac{2}{(p-1)}
\end{equation}
and so we have $\frac2p\le k_0\le 2$. But if $k_0=1$ then
$(p-1)k_1=1$, which is clearly a contradiction. Hence we should have
$k_0=2$, which implies that $k_1=k_2=\ldots=k_{p-1}=0$. This
completes the proof of Lemma \ref{lem3.2}.
\end{prf}

Now it is immediate to obtain the following corollary.

\begin{cor} \label{cor3.3}
Let $X$ be a homotopy $K3$ surface, and let $X$ admit a periodic
diffeomorphism $\tau$ of order $3$ satisfying $b_+(X/\tau)=3$.
Assume that the spin number $\text{\rm Spin}(\hat\tau, X)$ is
non-negative. Then we have
\[
k_0=2\ \text{and}\ k_1=k_2=0.
\]
\end{cor}

\begin{prf}
If the prime order $p$ equals 3 then the spin number $\text{\rm
Spin}(\hat\tau, X)$ is rational by Corollary \ref{cor3.2}. Thus the
equation \eqref{eq3.5} has rational coefficients. Thus all the
coefficients of \eqref{eq3.5} should vanish. Therefore, we have
$k_1=k_2$, as required.
\end{prf}

Next we deal with the case that the spin number is negative.

\begin{lemma} \label{lem3.3}
Let $X$ be a homotopy $K3$ surface, and let $X$ admit a periodic
diffeomorphism $\tau$ of odd prime order $p$ satisfying
$b_+(X/\tau)=3$. Assume that the spin number $\text{\rm
Spin}(\hat\tau, X)$ is both rational and negative. Then we have
\[
k_0 \le 0\ \text{and}\ k_1=k_2=\ldots=k_{p-1}=\frac{2-k_0}{p-1}\ge
1.
\]
Moreover, the sum $k_0+k_1+\ldots+k_{p-2}$ is equal to $2-k_1$ which
is non-negative.
\end{lemma}

\begin{prf}
From the relation \eqref{eq3.8}, note that
\begin{equation*}
0 > \text{Spin}(\hat\tau, X)= \frac{p k_0}{(p-1)}-\frac{2}{(p-1)}.
\end{equation*}
Thus we have $k_0< 2/p$, i.e., $k_0 \le 0$. The other relation
concerning $k_1$, $k_2, \cdots, k_{p-1}$ are immediate from the
identity $k_1=k_2=\cdots =k_{p-1}$, $k_0+k_1+\cdots+k_{p-1}=2$,
and Remark \ref{rmk3.1}.
\end{prf}

In particular, if $p=3$ then $k_0=2- 2k_1 \ge -2$. Thus under the
assumption that the spin number is rational and negative we have
either $k_0=0, k_1=k_2=1$ or $k_0=-2, k_1=k_2=2$. This will be used
in the proof of Proposition \ref{prop4.1}.

\section{Vanishing Theorems} 

The aim of this section is prove the main Theorem \ref{thm1.2}. We
do so by proving a vanishing theorem analogous to a theorem of F.
Fang in \cite{Fa}.

Before proving the Proposition \ref{prop4.1}, we first need to
consider the case that the spin number is both rational and
non-negative. In this case, it follows from Lemma \ref{lem3.2} that
$k_0=2$ and $k_1=k_2=\cdots=k_{p-1}=0$. But these data are exactly
what we can obtain for the trivial action on a homotopy $K3$
surface. Hence we cannot apply the arguments of this section to
obtain a new calculation of the Seiberg-Witten invariants, but only
the relation $\beta_f=a(1-t)^{m-1}$ for some integer $a$ holds. In
fact, if the arguments work in this case, we would have the
conclusion that there is no action of a cyclic group of odd prime
order $3$ on a homotopy $K3$ surface that acts trivially on $H^2_+$.
But there do exist examples of holomorphic actions of ${\bf Z}_3$ on
a projective $K3$ surface which is trivial on $H_+^2$, as we have
already seen in Section 1. On the other hand, the assumption that
the spin number is negative implies that the action is indeed
non-trivial, and the method of this section does not make any
contradiction to the known calculation of Morgan and Szab\' o for
homotopy $K3$ surfaces in \cite{M-S}.

Now we prove the following proposition which will be crucial
throughout this paper.

\begin{prop} \label{prop4.1}
Let $M$ be a smooth closed oriented spin 4-manifold with
$b_1(M)=0$ and $b_+(M)\ge 2$. Let $\tau$ generate a spin action of
odd prime order $p$ such that $b_+(M)=b_+(M/\tau)$. Assume that
the virtual dimension of the Seiberg-Witten moduli space is zero.
If $k_0\le l=\frac12(b_+(M) - 1)$ and $k_1=k_2=\ldots=k_{p-1}$
then we have
\begin{equation*}
\beta_f=a(1+\xi+\ldots+\xi^{p-1})(1-t)^{m_0-1}(1-t\xi)^{m_1}
\ldots(1-t\xi^{p-1})^{m_{p-1}},\ a\in {\Bbb Z}.
\end{equation*}
\end{prop}

\begin{remark} \label{rmk4.1}

\begin{itemize}

\item[(1)] By the relation $k_0+k_1\ldots+k_{p-1}=2$, the
assumption of the Proposition \ref{prop4.1} implies
$k_1=k_2=\ldots=k_{p-1}=\frac{2-k_0}{p-1}$.

\item[(2)] We assume in the statement of Proposition \ref{prop4.1} that the
virtual dimension $d$ of the Seiberg-Witten moduli space is zero,
since this is the only case we have for the trivial spin$^c$
structure on homotopy $K3$ surfaces. For the method how to deal with
the general case, you can see the proof of Theorem 1 in Section 4 of
\cite{Fa}. Indeed, for the case that the virtual dimension $d$ of
the Seiberg-Witten moduli space is not necessarily zero, we need to
replace $m$ by $m-d$ in the proof below and to consider the image of
$\beta_f$ in the truncated ring
\[
\frac{R(S^1\times {\Bbb Z}_p)}{(1-t)^{m_0-d}(1-t\xi)^{m_1}
\ldots(1-t\xi^{p-1})^{m_{p-1}}}.
\]
Then we see that the proof below can be repeated \emph{verbatim}
without any difficulty.
\end{itemize}
\end{remark}

\begin{prf}
For the sake of simplicity, we give a proof only for the case $p=3$.
This will not only greatly simplify the notational complications in
the proof, but also convey our idea more quickly. The other case is
completely similar. Refer to Section 5 of \cite{Fa} to see how to
deal with the general case in more detail.

Recall that our spin action is of even type. Let us denote by
$f=f_\Lambda$ the $S^1\times {\Bbb Z}_3$-equivariant map
\begin{equation*}
f: S^{(V_\lambda\oplus {\Bbb R})\otimes {\Bbb C}}\wedge
S(U_\lambda)\to S^{W'_\lambda\oplus V_\lambda\oplus {\Bbb R}}.
\end{equation*}
which is induced from the Seiberg-Witten equations as in
\eqref{eq2.1}.

Now, applying the Adams $\psi$-operation to the $K$-theoretic
degree $\beta\in K_{S^1\times {\Bbb Z}_3}(S(U_\Lambda))$ of $f$,
we get
\begin{equation}\label{eq4.1}
\begin{split}
\psi^q(\beta)=&q^l \beta\cdot
(1+t+\ldots+t^{q-1})^{n_0}\cdot(1+t\xi+\ldots+t^{q-1}\xi^{q-1})^{n_1}\\
&\cdot (1+t\xi^2+\ldots+t^{q-1}\xi^{2(q-1)})^{n_2}.
\end{split}
\end{equation}

We will need to use the following lemma in \cite{Fa}:

\begin{lemma} \label{lem4.1}
\[
K_{S^1\times {\Bbb Z}_3}(S(U_\lambda))=\frac{R(S^1\times {\Bbb
Z}_3)}{(1-t)^{m_0}(1-t\xi)^{m_1} (1-t\xi^2)^{m_2}},
\]
where $R(S^1\times {\Bbb Z}_3)$ denotes the representation ring of
$S^1\times {\Bbb Z}_3$.
\end{lemma}

Then we can show the following lemma:

\begin{lemma} \label{lem4.2}
There exists $\beta^{(1)}\in R(S^1\times {\Bbb
Z}_3)/(1-t\xi)^{m_1}(1-t\xi^2)^{m_2}$ such that
$\beta=\beta^{(1)}(1-t)^{m_0}$.
\end{lemma}

\begin{prf}
Let $\beta=\sum_i \left(\sum_{j=0}^2 a^i_j \xi^j\right) T^i$ for
$T=1-t$. Applying the identity \eqref{eq4.1} with $q=2$, we get
\begin{equation*}
\begin{split}
\psi^2(\beta)&=\sum_i\left(\sum_{j=0}^2 a^i_j \xi^{2j}\right)(2T-T^2)^i\\
&=2^l \sum_i\left(\sum_{j=0}^2 a^i_j \xi^{j}\right)T^i
(2-T)^{n_0}(1+\xi-T\xi)^{n_1}(1+\xi^2-T\xi^2)^{n_2}.
\end{split}
\end{equation*}

If we compare the coefficients of $T^i$, we get
\begin{equation} \label{eq4.2}
2^i\left(\sum_{j=0}^2 a^i_j
\xi^{2j}\right)=2^{l+n_0}(1+\xi)^{n_1}(1+\xi^2)^{n_2}\left(\sum_{j=0}^2
a^i_j \xi^{j}\right).
\end{equation}
Now, put $\xi=1$ in the equation \eqref{eq4.2}. Then we get
\begin{equation*}
2^i\left( \sum_{j=0}^2 a^i_j\right) =2^{l+n}\left( \sum_{j=0}^2
a^i_j\right),
\end{equation*}
where $n=n_0+n_1+n_2$. Thus if $i< l+n$, $\sum_{j=0}^2 a^i_j=0$.
Since the virtual dimension $d$ of the Seiberg-Witten moduli space
is zero and so $l+n=m-1=(m_0-1)+m_1+m_2\ge m_0$, we conclude that
$\sum_{j=0}^2 a^i_j=0$ for all $i\le m_0-1$. In fact, we can show
that  $a^i_0=a^i_1=a^i_2=0$ for all $i\le m_0-1$. To see it, we use
an argument of Fang in \cite{Fa}. Using the relation
$(1+\nu)(1+\nu^2)=1$ and the equation \eqref{eq4.2}, we obtain
\begin{equation*}
2^{2i}\prod_{k=1}^2\left( \sum_{j=0}^2 a^i_j
\nu^{2kj}\right)=2^{2(l+n_0)}\prod_{k=1}^2\left( \sum_{j=0}^2
a^i_j \nu^{kj}\right).
\end{equation*}
Since $\prod_{k=1}^2\left( \sum_{j=0}^2 a^i_j
\nu^{2kj}\right)=\prod_{k=1}^2\left( \sum_{j=0}^2 a^i_j
\nu^{kj}\right)$, there exist a $k$ such that $\sum_{j=0}^2 a^i_j
\nu^{kj}=0$, provided $i<l+n_0$. Thus we should have
$a^i_0=a^i_1=a^i_2=0$ for all $i\le m_0-1$.  Hence the image of
$\beta$ in $R(S^1\times {\Bbb Z}_3)/(1-t)^{m_0}$ is zero. Thus it
follows from Lemma \ref{lem4.1} and a simple argument similar to
Lemma 4.2 in \cite{Fa} that there exist a $\beta^{(1)}\in
R(S^1\times {\Bbb Z}_3)/(1-t\xi)^{m_1}(1-t\xi^2)^{m_2}$ such that
$\beta=\beta^{(1)}(1-t)^{m_0}$. This completes the proof.
\end{prf}

Similarly, we can show the following lemma whose proof is quite
similar to the case above.

\begin{lemma} \label{lem4.3}
There exists $\beta^{(2)}\in R(S^1\times {\Bbb
Z}_3)/(1-t\xi^2)^{m_2}$ such that
$\beta=\beta^{(2)}(1-t)^{m_0}(1-t\xi)^{m_1}$.
\end{lemma}

\begin{prf}
Let $\beta^{(1)}=\sum_i \left(\sum_{j=0}^2 b^i_j \xi^j\right) S^i$
for $S=1-t\xi$. Applying the identity \eqref{eq4.1} with $q=2$
again, it is easy to get
\begin{equation*}
\begin{split}
&(1+\xi^2-S\xi^2)^{k_0-1}\left(\sum_i \left( \sum_{j=0}^2 b^i_j
\xi^{2j}\right) (2S-S^2)^i\right)\\
&=2^l\left(\sum_i \left( \sum_{j=0}^2 b^i_j \xi^j\right)
S^i\right)(2-S)^{n_1}(1+\xi-S\xi)^{n_2}.
\end{split}
\end{equation*}

Now, comparing the coefficients of $S^i$, we get
\begin{equation} \label{eq4.3}
2^i(1+\xi^2)^{k_0-1}\left(\sum_{j=0}^2 b^i_j
\xi^{2j}\right)=2^{l+n_1}(1+\xi)^{n_2}\left(\sum_{j=0}^2 b^i_j
\xi^{j}\right).
\end{equation}
By putting $\xi=1$ in the equation \eqref{eq4.3}, we get
\begin{equation*}
2^{k_0+i-1}\left( \sum_{j=0}^2 b^i_j\right) =2^{l+n_1+n_2}\left(
\sum_{j=0}^2 b^i_j\right).
\end{equation*}
Thus if $k_0+i-1< l+n-n_0$, $\sum_{j=0}^2 b^i_j=0$. Since
$l+n=m-1=(m_0-1)+m_1+m_2$, the inequality $k_0+i-1< l+n-n_0$ is
equivalent to $i< m_1+m_2$. Thus we can conclude that
$\sum_{j=0}^2 b^i_j=0$ for all $i\le m_1-1$. In fact, an argument
similar to the previous case shows that $b^i_0=b^i_1=b^i_2$ for
all $i\le m_1-1$, and the rest of the proof is exactly same as
above. Thus we leave it to the reader. This completes the proof.
\end{prf}
\smallskip

Finally, let $\beta^{(2)}=\sum_i \left(\sum_{j=0}^2 c^i_j
\xi^j\right) Y^i$ for $Y=1-t\xi^2$. By the identity \eqref{eq4.1}
with $q=2$ again, we obtain
\begin{equation} \label{eq4.4}
\begin{split}
&(1+t)^{k_0-1}(1+t\xi)^{k_1}\left(\sum_i \left( \sum_{j=0}^2 c^i_j
\xi^{2j}\right) (1-t^2\xi)^i\right)\\
&=2^l\left(\sum_i \left( \sum_{j=0}^2 c^i_j \xi^j\right)
(1-t\xi^2)^i\right)(1+t^2\xi)^{n_2}.
\end{split}
\end{equation}
If we compare the coefficients of $Y^i$, we obtain
\begin{equation*}
(1+\xi)^{k_0-1}(1+\xi^2)^{k_1} 2^i \sum_{j=0}^2 c_j^i \xi^{2j} =2^l
2^{n_2} \sum_{j=0}^2 c_j^i \xi^j.
\end{equation*}
Since $k_0+k_1+i-1< l+n_2=l+n-n_0-n_1=(m-1)-n_0-n_1=m_2-1+k_0+k_1$,
for $i\le m_2-1$ we have $\sum_{j=0}^2 c_j^i=0$. On the other hand,
a similar argument as in Lemma \ref{lem4.2} shows that
\[
\sum_{j=0}^2 c_j^i \nu^{kj}=0
\]
for $i< l+n_2=m_2-1+k_0+k_1$ and some $k$. Since $k_0+k_1 \ge 0$
by Lemma \ref{lem3.3}, we have proved that $c^i_0=c^i_1=c^i_2=0$
for all $i\le m_2-2$.

Finally, if we substitute $Y=1-t\xi^2$ to the equation
\eqref{eq4.4}, we get
\begin{equation} \label{eq4.5}
\begin{split}
&(1+\xi-Y\xi)^{k_0-1}(1+\xi^2-Y\xi^2)^{k_1}(c^{m_2-1}_0
+c^{m_2-1}_1\xi^2+c^{m_2-1}_2\xi)(2Y-Y^2)^{m_2-1}\\
&=2^l
(c_0^{m_2-1}+c_1^{m_2-1}\xi+c_2^{m_2-1}\xi^2)Y^{m_2-1}(2-Y)^{n_2}.
\end{split}
\end{equation}
Now, comparing the coefficients of $Y^{m_2-1}$ in the equation
\eqref{eq4.5}, we get
\begin{equation} \label{eq4.6}
\begin{split}
&2^{m_2-1}(c^{m_2-1}_0+c^{m_2-1}_1\xi^2+c^{m_2-1}_2\xi)(1+\xi)^{k_0-1}(1+\xi^2)^{k_1}\\
&=2^{l+n_2}(c^{m_2-1}_0+c^{m_2-1}_1\xi+c^{m_2-1}_2\xi^2).
\end{split}
\end{equation}
Since $l+n_2-m_2+1=k_0+k_1$, it follows from \eqref{eq4.6} and the
relation $(1+\nu)(1+\nu^2)=1$ that we have
\begin{equation*}
\prod_{k=1}^2 \left(\sum_{j=0}^2 c^{m_2-1}_j \nu^{2kj} \right)=
2^{2(2-k_1)} \prod_{k=1}^2 \left(\sum_{j=0}^2 c^{m_2-1}_j
\nu^{kj}\right).
\end{equation*}
Thus if $k_1$ is not equal to $2$, then we have $\sum_{j=0}^2
c^{m_2}_j\nu^{kj}=0$ for some $k=1$ or $2$. Since $\nu$ is a
generic generator of the cyclic group ${\bf Z}_3$, this implies
$c^{m_2-1}_0=c^{m_2-1}_1=c^{m_2-1}_2$. Thus $\beta$ should be of
the form
\begin{equation} \label{eq4.7}
\beta=a(1+\xi+\xi^2)(1-t)^{m_0}(1-t\xi)^{m_1}(1-t\xi^2)^{m_2-1},\
a\in {\Bbb Z},
\end{equation}
as asserted. On the other hand, if $k_1$ is equal to $2$ then so
is $k_2$, and thus the spin number $\text{Spin}(\widehat{e^{2\pi
i/3} \tau}, X)$ is non-negative. But this implies that we should
have $k_0=k_2=0$ by a similar argument as in Lemma \ref{lem3.2},
which is a contradiction. This completes the proof of Proposition
\ref{prop4.1}.
\end{prf}

Now we are ready to prove our main theorem of this section.

\begin{theorem} \label{thm4.1}
Let $X$ be a homotopy $K3$ surface, and let $\tau:X\to X$ be a
periodic diffeomorphism $\tau$ of odd prime order $p$ satisfying
$b_+(X/\tau)=3$. Assume that the spin number $\text{\rm
Spin}(\hat\tau, X)$ is both rational and negative. Then the
Seiberg-Witten invariant for the trivial spin$^c$ structure vanishes
identically.
\end{theorem}

\begin{prf}
Since the difference of two spin structures is naturally an element
of $H^1(X, {\Bbb Z}_2)=0$, $\tau$ clearly generates a spin action of
odd prime order $p$. Moreover, since in our case we may assume that
$k_0\le 0$ and $k_1=k_2=\ldots=k_{p-1}$ by Lemma \ref{lem3.3}, all
the conditions in Proposition \ref{prop4.1} are satisfied. Thus we
can conclude from Proposition \ref{prop4.1} that $\beta_f$ should be
of the form
\begin{equation*}
\beta_f=a(1+\xi+\ldots+\xi^{p-1})(1-t)^{m_0}(1-t\xi)^{m_1}\ldots
(1-t\xi^{p-1})^{m_{p-1}-1},\ a\in {\Bbb Z}.
\end{equation*}
Now we claim that $\beta_f$ vanishes identically. To show this,
for our convenience we let
\begin{equation*}
\beta_f=\gamma(1-t)^{m_0}(1-t\xi)^{m_1}\ldots(1-t\xi^{p-1})^{m_{p-1}-1},\
a\in {\Bbb Z},
\end{equation*}
where $\gamma=a(1+\xi+\ldots+\xi^{p-1})$. By putting $\beta_f$
into the identity \eqref{eq4.1}, it is straightforward to see
\begin{equation} \label{eq4.8}
\psi^q(\gamma)(1+t+\ldots+t^{q-1})^{k_0}\cdots(1+t\xi^{p-1}+
\cdots+t^q\xi^{q(p-1)})^{k_{p-1}-1}=q^l \gamma.
\end{equation}
By plugging $\xi=1$ and $q=p$ (say) into \eqref{eq4.8} and using the
identity $t=1-T$ we should have
\begin{equation*}
ap(1+t\ldots+t^{p-1})=ap(p-\frac{(p-1)p}{2}T+\ldots)=p^2 a.
\end{equation*}
(It will not lose any generality that $m$ is sufficiently large, so
that $m>1$ in case that the virtual dimension $d$ of the
Seiberg-Witten moduli space is zero. e.g., see Section 4 of
\cite{Fa}.) Therefore, we see that $a$ should be zero. This
completes the proof.
\end{prf}

Now we are ready to prove one of the main results.

\begin{theorem} \label{thm4.2}
Let $X$ be a homotopy $K3$ surface, and let $\tau: X\to X$ be a
periodic diffeomorphism of odd prime order $p$. Assume that the
spin number $\text{\rm Spin}(\hat\tau, X)$ is rational and
negative. Then $\tau$ cannot act trivially on the self-dual part
$H^{2}_+(X; {\bf R})$ of the second cohomology group.
\end{theorem}

\begin{prf}
Suppose that $\tau$ acts trivially on the self-dual part
$H^{2}_+(X; {\bf R})$ of the second cohomology group. Then we have
$b_+(X/{\bf Z}_p)=3$. Since the spin number $\text{\rm
Spin}(\hat\tau, X)$ is rational and negative by the assumption, we
can apply Theorem \ref{thm4.1} to derive a contradiction to the
theorem of Morgan and Szab\' o in \cite{M-S}. This completes the
proof.
\end{prf}

\section{Applications} 

The aim of this section is to give just a few applications of our
main results and finally prove Theorem \ref{thm1.1}. To do so,
recall some necessary facts about psuedofree actions. An action is
called \emph{pseudofree} if it is free on the compliment of a
discrete subset. Now assume at the moment that the action of ${\bf
Z}_3$ is psuedofree. When we fix a generator of a cyclic group ${\bf
Z}_3$, the representation at each isolated fixed point can be
determined by a pair of non-zero integers $(\alpha, \beta)$ modulus
3 which is well-defined up to order and signs. Thus there are only
two types of $(\alpha, \beta)$; $(1,2)$ and $(1,1)$. Let $f_1$
(resp. $f_2$) be the number of fixed points of type $(1,2)$ (resp.
$(1,1)$). Then by the Atiyah-Singer $G$-signature theorem it is easy
to find
\begin{equation} \label{eq5.1}
3\sigma(X/{\bf Z}_3)= \sigma(X) + \frac{2}{3}(f_1- f_2).
\end{equation}
Since we have
\begin{equation} \label{eq5.2}
3 \chi(X/{\bf Z}_3)-2(f_1+f_2)=\chi(X),
\end{equation}
it follows from \eqref{eq5.1} and \eqref{eq5.2} that we have
\[
4f_1+2f_2 = 9b_+(X/{\bf Z}_3)-3.
\]
But $b_+(X/{\bf Z}_3)$ is either 1 or 3. Hence we have $2f_1+f_2=3$
or $2f_1+f_2=12$. On the other hand, notice that since the signature
$\sigma(X/{\bf Z}_3)$ is always integer, $f_1 -f_2\equiv 6$ mod 9.
This fact will not be used in this paper, but will be useful when we
try to classify cyclic group actions of order $3$ on $K3$ surfaces.

\begin{theorem} \label{thm5.1}
Let $X$ be a homotopy $K3$ surface, and let $\tau: X\to X$ be a
periodic diffeomorphism of order $3$. Assume that the fixed point
set is isolated. Then $\tau$ cannot act trivially on cohomology.
\end{theorem}

\begin{prf}
Suppose that $\tau$ acts trivially on cohomology. Then it follows
from \eqref{eq5.1} that $-16=\sigma(X)=\frac{1}{3}(f_1 -f_2)$. But
the spin number $\text{Spin}(\hat\tau, X)$ is the same as
$\frac{1}{3}(f_1 -f_2)$ by the Atiyah-Singer $G$-spin theorem.
Thus this number is always negative. Now if we combine the Theorem
\ref{thm4.2} together with Corollary \ref{cor3.2}, we are done.
This completes the proof.
\end{prf}

Moreover, we can show that a homotopy $K3$ surface admits a
periodic diffeomorphism of order $3$ which acts trivially on
cohomology only if its spin number is negative. Hence we can give
the following

\begin{theorem} \label{thm5.2}
Let $X$ be a homotopy $K3$ surface, and let $\tau: X\to X$ be a
periodic diffeomorphism of order $3$. Then $\tau$ cannot act
trivially on cohomology.
\end{theorem}

\begin{prf}
By the Atiyah-Singer $G$-signature theorem (e.g., see Proposition
6.18 of \cite{A-S}) we have
\begin{equation} \label{eq5.3}
3\sigma(X/{\bf Z}_3)=\sigma(X)+\sum_k \sum_{l=1}^2 \csc^2(\pi l/3)
\langle [F_k], [F_k] \rangle + \frac{2}{3} (f_1 - f_2).
\end{equation}
Moreover, since the spin number $\text{Spin}(\hat\tau, X)$ is
given by
\begin{equation} \label{eq5.4}
\frac{1}{4} \sum_k \cos({\pi /3}) \csc^2({\pi /3}) \langle [F_k],
[F_k] \rangle + \frac{1}{3} (f_1 -f_2)
\end{equation}
by the Atiyah-Singer $G$-spin theorem, we have
\begin{equation} \label{eq5.5}
\begin{split}
3 \text{Ind}_{{\bf Z}_3} D &= -\frac{\sigma(X)}{8} + 2
\text{Spin}(\hat\tau, X)\\
&=-\frac{\sigma(X)}{8}+\frac{1}{4} \sum_k \sum_{l=1}^2 (-1)^l
\cos({\pi l/3}) \csc^2({\pi l/3}) \langle [F_k], [F_k] \rangle\\
&+ \frac{2}{3} (f_1 -f_2).
\end{split}
\end{equation}
Here the signs in the spin number are determined as in the paper
\cite{AH} of Atiyah and Hirzebruch.

Now assume that $\tau$ acts trivially on cohomology. Then it
follows from \eqref{eq5.3} and \eqref{eq5.5} that we have
\begin{equation} \label{eq5.6}
k_0= 2+ \frac{1}{4}(f_1 -f_2),
\end{equation}
where we use the identity $k_0=\text{Ind}_{{\bf Z}_3} D$ and
\[
\sum_{l=1}^{p-1} (-1)^l \cos({\pi l/p}) \csc^2({\pi l /p}) =
-\frac12 \sum_{l=1}^{p-1} \csc^2({\pi l/p}).
\]
If $f_1=f_2$ then it follows from \eqref{eq5.4} that the spin number
$\text{Spin}(\hat\tau, X)$ becomes
\[
\frac16 \sum_k \langle [F_k], [F_k] \rangle.
\]
It is known in \cite{Ed2} and \cite{Mc} that the fixed point set of
each group element consists of 2-spheres and/or isolated fixed
points. Hence by the adjunction inequality, the self-intersection
number $\langle [F_k], [F_k]\rangle$ is always less than or equal to
0. Thus the spin number $\text{Spin}(\hat\tau, X)$ is non-positive
in this case. But if $f_1=f_2$, then it also follows from
\eqref{eq5.6} that $k_0=2$, which implies that the spin number is
positive. Hence we can conclude that this case does not happen.
Notice that $\text{Spin}(\hat\tau, X)$ cannot be zero (e.g., see
\eqref{eq3.8}).

Next if $f_1\ne f_2$, then using $k_0\le 2$ by Lemma \ref{lem3.1} we
have $k_0 <2$, which implies that the spin number
$\text{Spin}(\hat\tau, X)$ is negative. Hence in this case we can
apply Theorem \ref{thm4.1} to obtain the vanishing of the
Seiberg-Witten invariant for the trivial spin$^c$ structure. This is
a contradiction, which completes the proof.
\end{prf}

In fact, there exists an example of a $K3$ surface admitting a
homologically nontrivial action of ${\bf Z}_3$ whose spin number
is positive.

\begin{example} Consider the Fermat quartic surface $X$ which is
defined by the equation $\sum_{1}^4 z_i^4 =0$ in ${\bf CP}^3$.
Then define an action of ${\bf Z}_3$ on $X$ generated by
\[
[z_1, z_2, z_3, z_4]\mapsto [z_1, z_3, z_4, z_2].
\]
Then the ${\bf Z}_3$-action has 6 isolated fixed points: four
points of the form $[1,a,a,a]$ with $1+3a^4=0$ and two more points
of the form $[0,1,b^2, b]$ with $b^2+b+1=0$. In our case
$f_1+f_2=6$. Thus using the identity $2f_1+f_2=12$, we have
$f_1=6$ and $f_2=0$. On the other hand, the spin number
$\text{spin}(\hat\tau, X)$ is given by $\frac{1}{3}(f_1 -f_2)$, if
the fixed point set is isolated. Thus in our case the spin number
$\text{Spin}(\hat\tau, X)$ is $2$ which is positive. It is also
easy to see that $b_+(X/{\bf Z}_3)=3$ and $b_-(X/{\bf Z}_3)=7$. So
the group action is nontrivial on $H^2(X, {\bf R})$.
\end{example}


\end{document}